\documentclass{gtart}

\def\ifplaintex{\expandafter\ifx\csname documentclass\endcsname\relax}

\ifplaintex 
\else
\fi

\expandafter\ifx\csname beginpicture\endcsname\relax
\expandafter\ifx\csname documentclass\endcsname\relax
\input pictex \else
\input prepictex \input pictex \input postpictex \fi\fi

\def\gt{{\mathsurround=0pt\it $\cal G\mskip-2mu$eometry \&\ 
$\cal T\!\!$opology}}        %  journal title in recommended style

\def\gtp{{\mathsurround=0pt\it $\cal G\mskip-2mu$eometry \&\ 
$\cal T\!\!$opology $\cal P\!$ublications}}  % GT publications

\def\lognumber#1{\def\thelognumber{#1}}
\def\volumenumber#1{\def\thevolumenumber{#1}}
\def\papernumber#1{\def\thepapernumber{#1}}
\def\volumeyear#1{\def\thevolumeyear{#1}}

\def\pagenumbers#1#2{\def\startpage{#1}\def\finishpage{#2}}
\def\published#1{\def\publishdate{#1}}
\def\proposed#1{\def\theproposer{#1}}
\def\seconded#1{\def\theseconders{#1}}
\def\received#1{\def\receiveddate{#1}}
\def\revised#1{\def\reviseddate{#1}}
\def\accepted#1{\def\accepteddate{#1}}

\long\def\asciiabstract#1{\long\def\theasciiabstract{#1}}

\let\\\par\let\thelognumber\relax
\let\thevolumenumber\relax\let\thepapernumber\relax
\let\thevolumeyear\relax\let\thesamplenumber\relax\let\startpage\relax
\let\finishpage\relax\let\publishdate\relax\let\receiveddate\relax
\let\reviseddate\relax\let\accepteddate\relax\let\theasciititle\relax
\let\theasciiauthors\relax
\let\theasciiabstract\relax
\let\theasciiemail\relax\let\theshortauthors\relax\let\theshorttitle\relax

\long\def\maketitlep{   % start of definition of \maketitlep

\count0=\startpage

\gt\hfill      %   Journal title (top left) 
\beginpicture
\setcoordinatesystem units <0.33truein, 0.33truein> point at 2.2 0.9
\setplotsymbol ({$\cal G$})
\plotsymbolspacing=9truept
\circulararc 315 degrees from 0 1 center at 0 0
\setplotsymbol ({$\cal T$})
\circulararc 315 degrees from 1 -1 center at 1 0
\endpicture
\break
{\small\ifx\thesamplenumber\relax % sample?  
Volume \else Sample
\fi\thevolumenumber\ (\thevolumeyear)
\startpage--\finishpage\nl
Published: \publishdate}
\vglue 0.5truein plus 0.4fil minus 0.1truein

{\parskip=0pt\leftskip 0pt plus 1fil\def\\{\par\smallskip}{\ifplaintex\large
\else\Large\fi\bf\thetitle}\par\medskip}   

\vglue 0pt plus 0.1fil 

{\parskip=0pt\leftskip 0pt plus 1fil\def\\{\par}{\sc\theauthors}
\par\medskip}

\vglue 0pt plus 0.1fil 

{\small\parskip=0pt\let\newline\\
{\leftskip 0pt plus 1fil\def\\{\par}{\sl\theaddress}\par}
\expandafter\ifx\theemail\relax    % email address?
\relax\else\vglue 5pt plus 0.02fil minus 2pt\def\\{\stdspace{\rm 
and}\stdspace} 
\cl{Email:\stdspace\tt\theemail}\fi
\ifx\theurl\relax                  % URL given?
\relax\else\vglue 5pt plus 0.02fil minus 2pt\def\\{\stdspace{\rm 
and}\stdspace}
\cl{URL:\stdspace\tt\theurl}\fi\par}

\vglue 7pt plus 0.3fil minus 3pt

{\bf Abstract}
\vglue 5pt plus 0.1fil minus 2pt

\theabstract

\vglue 7pt plus 0.3fil minus 3pt

{\bf AMS Classification numbers}\quad Primary:\quad \theprimaryclass

Secondary:\quad \thesecondaryclass

\vglue 5pt plus 0.3fil minus 2pt

{\bf Keywords}\quad \thekeywords

\vglue 10pt plus 0.5fil minus 5pt

{\small  Proposed: \theproposer\hfill Received: \receiveddate\nl
Seconded: \theseconders\hfill 
\ifx\reviseddate\relax                         % paper revised?
Accepted: \accepteddate                        % no
\else
Revised: \reviseddate                          % yes
\fi}
\eject
}       %  end of definition of \maketitlep

\let\maketitlepage\maketitlep
\let\maketitle\maketitlepage

\font\phead=cmsl9 scaled 950
\font=cmsl9 scaled 1050
\font\pnum=cmbx10 scaled 913
\font\lnum=cmbx10 
\font\pfoot=cmsl9 scaled 950
\font=cmsl9 scaled 1050
\ifplaintex
\headline{\vbox to 0pt{\vskip -4.5mm\line{\small\phead\ifnum
\count0=\startpage ISSN 1364-0380 (on line)
1465-3060 (printed) \hfill {\pnum\folio}\else\ifodd\count0\def\\{ }% 
\ifx\theshorttitle\relax\thetitle\else\theshorttitle\fi\hfill{\pnum\folio}
\else\def\\{ and }{\pnum\folio}\hfill\ifx\theshortauthors\relax\theauthors
\else\theshortauthors\fi\fi\fi}\vss}}
\footline{\vbox to 0pt{\vglue 0mm\line{\small\pfoot\ifnum\count0=\startpage
\copyright\ \gtp\hfill\else
\gt, Volume \thevolumenumber\ (\thevolumeyear)\hfill\fi}\vss
}}
\else
\makeatletter
\def\@oddhead{{\small\lhead\ifnum\count0=\startpage ISSN 1364-0380 (on line)
1465-3060 (printed) \hfill {\lnum\number\count0}\else\ifodd\count0
\def\\{ }\ifx\theshorttitle\relax \thetitle \else\theshorttitle\fi\hfill
{\lnum\number\count0}\else\def\\{ and }{\lnum\number\count0}
\hfill\ifx\theshortauthors\relax 
\theauthors\else\theshortauthors\fi\fi\fi}}\def\@evenhead{\@oddhead}
\def\@oddfoot{\small\lfoot\ifnum\count0=\startpage\copyright\ \gtp\hfill\else
\gt, Volume \thevolumenumber\ (\thevolumeyear)\hfill\fi}
\def\@evenfoot{\@oddfoot}
\makeatother
\fi

\newwrite\gtoutfile
\long\gdef\makeheadfile{  %%% start of definition of \makeheadfile
{\def\\{, }\def\s{ }
\immediate\openout\gtoutfile head.xxx
\immediate\write\gtoutfile{To: math@arxiv.org}
\immediate\write\gtoutfile{Subject: put or rep NNNNN:pppp}
\immediate\write\gtoutfile{--text follows this line--}
\immediate\write\gtoutfile{Proxy-for: \ifx\theasciiauthors\relax
\theauthors\else\theasciiauthors\fi\s<\ifx\theasciiemail\relax\theemail\else\theasciiemail\fi>}
\immediate\write\gtoutfile{\noexpand\\}
\immediate\write\gtoutfile{Authors: \ifx\theasciiauthors\relax
\theauthors\else\theasciiauthors\fi}
{\def\\{ }\immediate\write\gtoutfile{Title: \ifx\theasciititle\relax
\thetitle\else\theasciititle\fi}}
\immediate\write\gtoutfile{Subj-class: GT or SG or MG etc}
\immediate\write\gtoutfile{MSC-class: \theprimaryclass\ifx\thesecondaryclass\relax\else, \thesecondaryclass\fi}
\immediate\write\gtoutfile{Journal-ref: Geom. Topol. \thevolumenumber
(\thevolumeyear) \startpage-\finishpage}
\immediate\write\gtoutfile{Comments: Published by Geometry and Topology at}
\immediate\write\gtoutfile{\s\s http://www.maths.warwick.ac.uk/gt/GTVol\thevolumenumber/paper\thepapernumber.abs.html}
\immediate\write\gtoutfile{\noexpand\\}
\immediate\write\gtoutfile{}
\ifx\theasciiabstract\relax
\immediate\write\gtoutfile{\theabstract}\else
\immediate\write\gtoutfile{\theasciiabstract}\fi
\immediate\write\gtoutfile{}
\immediate\write\gtoutfile{\noexpand\\}
\immediate\write\gtoutfile{}
\immediate\closeout\gtoutfile}}  %%% end of definition of \makeheadfile

\def\maketitlepage{\maketitlep\makeheadfile}
\let\maketitle\maketitlepage

\def\ifplaintex{\expandafter\ifx\csname documentclass\endcsname\relax}

\ifplaintex 
\else
\fi

\expandafter\ifx\csname beginpicture\endcsname\relax
\expandafter\ifx\csname documentclass\endcsname\relax
\input pictex \else
\input prepictex \input pictex \input postpictex \fi\fi

\def\gt{{\mathsurround=0pt\it $\cal G\mskip-2mu$eometry \&\ 
$\cal T\!\!$opology}}        %  journal title in recommended style

\def\gtp{{\mathsurround=0pt\it $\cal G\mskip-2mu$eometry \&\ 
$\cal T\!\!$opology $\cal P\!$ublications}}  % GT publications

\def\lognumber#1{\def\thelognumber{#1}}
\def\volumenumber#1{\def\thevolumenumber{#1}}
\def\papernumber#1{\def\thepapernumber{#1}}
\def\volumeyear#1{\def\thevolumeyear{#1}}

\def\pagenumbers#1#2{\def\startpage{#1}\def\finishpage{#2}}
\def\published#1{\def\publishdate{#1}}
\def\proposed#1{\def\theproposer{#1}}
\def\seconded#1{\def\theseconders{#1}}
\def\received#1{\def\receiveddate{#1}}
\def\revised#1{\def\reviseddate{#1}}
\def\accepted#1{\def\accepteddate{#1}}

\long\def\asciiabstract#1{\long\def\theasciiabstract{#1}}

\let\\\par\let\thelognumber\relax
\let\thevolumenumber\relax\let\thepapernumber\relax
\let\thevolumeyear\relax\let\thesamplenumber\relax\let\startpage\relax
\let\finishpage\relax\let\publishdate\relax\let\receiveddate\relax
\let\reviseddate\relax\let\accepteddate\relax\let\theasciititle\relax
\let\theasciiauthors\relax
\let\theasciiabstract\relax
\let\theasciiemail\relax\let\theshortauthors\relax\let\theshorttitle\relax

\long\def\maketitlep{   % start of definition of \maketitlep

\count0=\startpage

\gt\hfill      %   Journal title (top left) 
\beginpicture
\setcoordinatesystem units <0.33truein, 0.33truein> point at 2.2 0.9
\setplotsymbol ({$\cal G$})
\plotsymbolspacing=9truept
\circulararc 315 degrees from 0 1 center at 0 0
\setplotsymbol ({$\cal T$})
\circulararc 315 degrees from 1 -1 center at 1 0
\endpicture
\break
{\small\ifx\thesamplenumber\relax % sample?  
Volume \else Sample
\fi\thevolumenumber\ (\thevolumeyear)
\startpage--\finishpage\nl
Published: \publishdate}
\vglue 0.5truein plus 0.4fil minus 0.1truein

{\parskip=0pt\leftskip 0pt plus 1fil\def\\{\par\smallskip}{\ifplaintex\large
\else\Large\fi\bf\thetitle}\par\medskip}   

\vglue 0pt plus 0.1fil 

{\parskip=0pt\leftskip 0pt plus 1fil\def\\{\par}{\sc\theauthors}
\par\medskip}

\vglue 0pt plus 0.1fil 

{\small\parskip=0pt\let\newline\\
{\leftskip 0pt plus 1fil\def\\{\par}{\sl\theaddress}\par}
\expandafter\ifx\theemail\relax    % email address?
\relax\else\vglue 5pt plus 0.02fil minus 2pt\def\\{\stdspace{\rm 
and}\stdspace} 
\cl{Email:\stdspace\tt\theemail}\fi
\ifx\theurl\relax                  % URL given?
\relax\else\vglue 5pt plus 0.02fil minus 2pt\def\\{\stdspace{\rm 
and}\stdspace}
\cl{URL:\stdspace\tt\theurl}\fi\par}

\vglue 7pt plus 0.3fil minus 3pt

{\bf Abstract}
\vglue 5pt plus 0.1fil minus 2pt

\theabstract

\vglue 7pt plus 0.3fil minus 3pt

{\bf AMS Classification numbers}\quad Primary:\quad \theprimaryclass

Secondary:\quad \thesecondaryclass

\vglue 5pt plus 0.3fil minus 2pt

{\bf Keywords}\quad \thekeywords

\vglue 10pt plus 0.5fil minus 5pt

{\small  Proposed: \theproposer\hfill Received: \receiveddate\nl
Seconded: \theseconders\hfill 
\ifx\reviseddate\relax                         % paper revised?
Accepted: \accepteddate                        % no
\else
Revised: \reviseddate                          % yes
\fi}
\eject
}       %  end of definition of \maketitlep

\let\maketitlepage\maketitlep
\let\maketitle\maketitlepage

\font\phead=cmsl9 scaled 950
\font=cmsl9 scaled 1050
\font\pnum=cmbx10 scaled 913
\font\lnum=cmbx10 
\font\pfoot=cmsl9 scaled 950
\font=cmsl9 scaled 1050
\ifplaintex
\headline{\vbox to 0pt{\vskip -4.5mm\line{\small\phead\ifnum
\count0=\startpage ISSN 1364-0380 (on line)
1465-3060 (printed) \hfill {\pnum\folio}\else\ifodd\count0\def\\{ }% 
\ifx\theshorttitle\relax\thetitle\else\theshorttitle\fi\hfill{\pnum\folio}
\else\def\\{ and }{\pnum\folio}\hfill\ifx\theshortauthors\relax\theauthors
\else\theshortauthors\fi\fi\fi}\vss}}
\footline{\vbox to 0pt{\vglue 0mm\line{\small\pfoot\ifnum\count0=\startpage
\copyright\ \gtp\hfill\else
\gt, Volume \thevolumenumber\ (\thevolumeyear)\hfill\fi}\vss
}}
\else
\makeatletter
\def\@oddhead{{\small\lhead\ifnum\count0=\startpage ISSN 1364-0380 (on line)
1465-3060 (printed) \hfill {\lnum\number\count0}\else\ifodd\count0
\def\\{ }\ifx\theshorttitle\relax \thetitle \else\theshorttitle\fi\hfill
{\lnum\number\count0}\else\def\\{ and }{\lnum\number\count0}
\hfill\ifx\theshortauthors\relax 
\theauthors\else\theshortauthors\fi\fi\fi}}\def\@evenhead{\@oddhead}
\def\@oddfoot{\small\lfoot\ifnum\count0=\startpage\copyright\ \gtp\hfill\else
\gt, Volume \thevolumenumber\ (\thevolumeyear)\hfill\fi}
\def\@evenfoot{\@oddfoot}
\makeatother
\fi

\newwrite\gtoutfile
\long\gdef\makeheadfile{  %%% start of definition of \makeheadfile
{\def\\{, }\def\s{ }
\immediate\openout\gtoutfile head.xxx
\immediate\write\gtoutfile{To: math@arxiv.org}
\immediate\write\gtoutfile{Subject: put or rep NNNNN:pppp}
\immediate\write\gtoutfile{--text follows this line--}
\immediate\write\gtoutfile{Proxy-for: \ifx\theasciiauthors\relax
\theauthors\else\theasciiauthors\fi\s<\ifx\theasciiemail\relax\theemail\else\theasciiemail\fi>}
\immediate\write\gtoutfile{\noexpand\\}
\immediate\write\gtoutfile{Authors: \ifx\theasciiauthors\relax
\theauthors\else\theasciiauthors\fi}
{\def\\{ }\immediate\write\gtoutfile{Title: \ifx\theasciititle\relax
\thetitle\else\theasciititle\fi}}
\immediate\write\gtoutfile{Subj-class: GT or SG or MG etc}
\immediate\write\gtoutfile{MSC-class: \theprimaryclass\ifx\thesecondaryclass\relax\else, \thesecondaryclass\fi}
\immediate\write\gtoutfile{Journal-ref: Geom. Topol. \thevolumenumber
(\thevolumeyear) \startpage-\finishpage}
\immediate\write\gtoutfile{Comments: Published by Geometry and Topology at}
\immediate\write\gtoutfile{\s\s http://www.maths.warwick.ac.uk/gt/GTVol\thevolumenumber/paper\thepapernumber.abs.html}
\immediate\write\gtoutfile{\noexpand\\}
\immediate\write\gtoutfile{}
\ifx\theasciiabstract\relax
\immediate\write\gtoutfile{\theabstract}\else
\immediate\write\gtoutfile{\theasciiabstract}\fi
\immediate\write\gtoutfile{}
\immediate\write\gtoutfile{\noexpand\\}
\immediate\write\gtoutfile{}
\immediate\closeout\gtoutfile}}  %%% end of definition of \makeheadfile

\def\maketitlepage{\maketitlep\makeheadfile}
\let\maketitle\maketitlepage

\lognumber{202}

\volumenumber{6}
\papernumber{14}
\volumeyear{2002}
\pagenumbers{403}{408}
\received{21 August 2001}
\revised{21 April 2002}
\accepted{22 August 2002}
\published{5 September 2002}
\proposed{Cameron Gordon}

\seconded{Ronald Stern, Walter Neumann}

\def\zz{{\bf Z}}
\def\qq{{\bf Q}}
  
\def\calc{\mathcal{C}}
\def\calg{\mathcal{G}}

\newtheorem{theorem}{Theorem}[section]
\newtheorem{lemma}[theorem]{Lemma}
\newtheorem{corollary}[theorem]{Corollary}

\theoremstyle{definition}

\theoremstyle{remark}

\begin{document}

\title{Seifert forms and concordance}
\author{Charles Livingston}
\address{Department of Mathematics, Indiana University\\Bloomington, 
IN 47405, USA}
\email{livingst@indiana.edu}

\begin{abstract}
If a knot $K$ has Seifert matrix $V_K$ and has a prime power cyclic
branched cover that is not a homology sphere, then there is an
infinite family of non--concordant knots having Seifert matrix $V_K$.
\end{abstract}
\asciiabstract{
If a knot K has Seifert matrix V_K and has a prime power cyclic
branched cover that is not a homology sphere, then there is an
infinite family of non-concordant knots having Seifert matrix V_K.}

\primaryclass{57M25}
\secondaryclass{57N70} 

\keywords {Concordance, Seifert matrix, Alexander polynomial}

\maketitlepage

%%%%%%%%%%%%%%%%%%%%%%%%%%%%%%%SECTION%%%%%%%%%%%%%%%%%%%%%%%%%%%%%%%

\section{Introduction}

Levine's homomorphism  $\psi\co\calc \to \calg$ from the concordance group
of knots in
$S^3$ to the algebraic concordance group of Seifert matrices (defined in \cite{le1}) has an
infinitely generated kernel, as proved by Jiang
\cite{ji}.   It follows that every algebraic concordance
class can be represented by an infinite family of non--concordant knots.  However, it is also the case
that every class in
$\calg$ can be represented by an infinite number of distinct Seifert matrices,  so Jiang's result
alone tells us nothing about whether a given Seifert matrix can arise from non--concordant knots.  In
fact, all the knots in the kernel of
$\phi$ identified by Jiang have distinct Seifert forms.  

Examples of non--slice, algebraically slice, knots quickly yield    pairs of non--concordant knots
with the same Seifert matrix. Beyond this  nothing has been known regarding the extent to
which the Seifert matrix of a knot might determine its concordance class.  We prove the following.

\begin{theorem}\label{maintheorem} If a knot $K$ has Seifert matrix $V_K$ and its
Alexander polynomial $\Delta_K(t)$  has an irreducible factor that is not a cyclotomic polynomial
$\phi_n$ with
$n$ divisible by three distinct primes, then there is an infinite family
$\{K_i\}$ of non--concordant knots such that each
$K_i$ has Seifert matrix
$V_K$.
\end{theorem}

The condition on the Alexander polynomial  seems somewhat technical; we note three relevant facts.
First,  if the Alexander polynomial of the knot is trivial,  $\Delta_K(t) = 1$, then $K$ is
topologically slice
\cite{fr, fq}. Second we have:

\begin{theorem}\label{polytheorem}  All prime power cyclic branched covers of a knot $K$ are homology
spheres if and only if all nontrivial irreducible factors of $\Delta_K(t)$ are cyclotomic polynomials
$\phi_{n}(t)$ with
$n$ divisible by three distinct primes. All branched covers of $K$ are homology spheres if and only if
$\Delta_K(t) = 1$. \end{theorem}

\noindent  Finally, we note that Taehee Kim \cite{k} has applied the recent advances in concordance
theory of \cite{cot} to prove that for each $n$ divisible by three distinct primes there is a knot
with
 $ \Delta_K(t) = \left(\phi_{n}(t)\right)^2$   for which there is an infinite family of
non--concordant knots having the same Seifert matrix. 

A good   reference  for the basic knot theory in this paper is \cite{rol}, for the
algebraic concordance group \cite{le1, le2} are the main references, and for Casson--Gordon invariants
references are \cite{cg1, gl1}. 

\medskip
{\bf Remark}\qua  
We have chosen to use Seifert matrices instead of Seifert forms to be
consistent with references \cite{le1,le2}.  A basis free approach
using Seifert forms could be carried out identically.

\section{Proof of Theorem \ref{maintheorem}}

Unless indicated, all homology groups are taken with integer coefficients.  

Let $F$ be a Seifert surface for $K$ with associated Seifert matrix $V_K$.  View $F$ as a disk with
$2g$ bands added and let $\{S_m\}_{ m = 1, \ldots , 2g}$, be a collection of unknotted circles, one
linking each of the bands. Let $K_i$ be the knot formed by replacing a tubular neighborhood of each
$S_m$ with a copy of the complement of  a knot $J_i$, identifying the meridian and longitude of $J_i$
with the longitude and meridian of the $S_m$, respectively.  The correct choice of the $J_i$ will
be identified in the proof. Replacing the $S_m$ with the knot complements has the effect of adding a
local knot to each band of
$F$.  The Seifert form of
$K_i$ is independent of the choice of
$J_i$.  Applying Theorem
\ref{polytheorem}, proved in the next section, we assume  the $p^k$--fold cyclic branched cover of
$S^3$ branched over $K$ has nontrivial homology. Denote this cover by $M(K)$ and let $q$ be a maximal
prime power divisor of $|H_1(M(K))|$.  

According to Casson and Gordon \cite{cg1}, if $K_i \# - K_j$ is slice (that is, if $K_i$ and $K_j$
are concordant) then for some nontrivial $\zz_q$--valued character $\chi$ on $H_1(M( K_i \# - K_j))$
the Casson--Gordon invariant $\sigma_1(\tau(K_i \# - K_j,\chi)) = 0$.  Using  the additivity of
Casson--Gordon invariants  (proved by Gilmer \cite{gi}), this equality can be rewritten as 
$\sigma_1(\tau(K_i ,\chi_i)) = \sigma_1(\tau(K_j ,\chi_j))$ where $\chi_i$ and $\chi_j$ are the
restrictions of $\chi$ to $H_1(M(K_i))$ and $H_1(M(K_j))$, respectively.  Notice that at least one of
$\chi_i$ and $\chi_j$ is nontrivial. Furthermore, since according to  \cite{cg1}  (see also
\cite{gi}) the set of characters for which the Casson--Gordon invariants must vanish is a metabolizer
for the
  linking form on $H^1(M(K_i \# - K_j), \qq/\zz)$, there are  such characters for which
$\chi_j$ must be nontrivial.  (If the metabolizer was contained in $H^1( M(K_i), \qq/\zz)$ then order
considerations would show that it equalled this summand, contradicting nonsingularity.)

Litherland's analysis \cite{lit2}  of companionship and Casson--Gordon invariants applies
directly to the case of knotting the bands in the Seifert surface  (see
also
\cite{gl1})).  Roughly stated, there is a correspondence between characters on $H_1(M(K))$ and on
$H_1(M(K_i))$; it then follows that the difference of the corresponding Casson--Gordon invariants is
determined by
$q$--signatures of
$J_i$:
$\sigma_{a/q}(J_i) = \mbox{sign}\left( (1-\omega) V_{J_i} + (1 -
\overline{\omega})V_{J_i}^t \right)$ where $\omega = e^{2\pi ai/q }$.  More precisely, it follows
readily from the results of \cite{lit2} and iteration that the equality of Casson--Gordon invariants
for
$K_i$ and
$K_j$ is given by 
$$(*)\ \ \ \ \ \  \sigma_1(\tau(K,\chi_i)) + \sum_l \sigma_{a_l/q}(J_i) =  \sigma_1(\tau(K,\chi_j)) +
\sum_l
\sigma_{b_l/q}(J_j).
$$

The two summations that appear have $2gp^k$ terms in them.  The values of the $a_l$ are given by the
values of $\chi_i$ on the $2gp^k$ lifts of the circles $S_m$ to $M(K)$.  Similar statements hold for
the $b_l$ and
$\chi_j$. Observe also that since the lifts of the $S_m$ generate $H_1(M(K))$ (see for instance
\cite{rol}) and at least one of
$\chi_i$ or
$\chi_j$ is nontrivial, at least one of the $a_l$ or $b_l$ is nontrivial.

A prime power branched cover of a knot is a rational homology sphere and hence $H_1(M(K))$ is
finite.  A short proof of this is given in the next section.  Hence, there is only a finite set of
characters to consider and 
$
\sigma_1(\tau(K,\chi_1))
$ lies in a bounded range, say $[-N_0,N_0]$.  If we can choose $J_1$ so that $\sum_l
\sigma_{a_l/q}(J_1)$  lies in a range  $[2N_0 +1,
N_1]$ (for some $N_1$ and for all possible sums with some $a_l \ne 0 \in \zz_q$) then it would follow
that
$K$ and
$K_1$ are not concordant. Similarly, by selecting each $J_{i+1}$ so that the sum lies in the range
$[2N_0 + N_i +1, N_{i+1}]$ we will have that the equality $(*)$ cannot hold for any pair $i$ and $j$
and the theorem is proved.

The desired $J_i$ are constructed by taking ever larger multiples of a knot $T$ for which
$\sigma_{a/q}(T) \ge 2$ for all $a \ne 0 \in \zz_q$.  Such a knot is given in the following lemma,
which completes the proof of Theorem \ref{maintheorem}.

\begin{lemma}  The $(2,q)$--torus knot $T_{2,q}$ has  $\sigma_{a/q}(T) \ge 2$ for all  $a \ne 0 \in
\zz_q$. \end{lemma}

\begin{proof}  The signature function of a knot $K$, $\mbox{sign}\left( (1-\omega) V_{K} + (1 -
\overline{\omega})V_{K} \right)$, has jumps only at roots of the Alexander polynomial, and if these
roots are simple the jump is either $\pm 2$ \cite{ma}.  The $(2,q)$--torus knot has cyclotomic
Alexander polynomial
$\phi_{2q}$ with $(q-1)/2$ simple roots on the upper unit circle in the complex plane.  Hence the
signature $\sigma_{-1}(T_{2,q} )\le q-1$.  On the other hand, this $-1$ signature is easily computed
from the standard rank $q-1$ Seifert form for $T_{2,q}$ to be exactly $q-1$, and so all the jumps must
be positive 2.  The first of these jumps occurs at a primitive $2q$--root of unity, so  all
$q$--signatures must be positive as desired.
\end{proof}

%%%%%%%%%%%%%%%%%%%%%%%%%%%%%%%SECTION%%%%%%%%%%%%%%%%%%%%%%%%%%%%%%%

\section{Proof of Theorem \ref{polytheorem}}

We have the following result of Fox \cite{fox0}  and include as a corollary a result used above.

\begin{theorem}  \label{polyproduct} If $M(K)$ is the $r$--fold cyclic branched cover of $S^3$
branched over $K$, then $$| H_1(M(K))|  =
\prod_{i=0 }^{ r-1}
\Delta_K(\zeta_r^i)$$ where $\zeta_r$ is a primitive $r$--root of unity. If the product is 0 then
$H_1(M(K))$ is infinite.
\end{theorem}

\begin{corollary}\label{rat} If $r$ is a prime power, then $M(K)$ is a rational homology sphere:
$ H_1(M(K), \qq)  = 0$.
\end{corollary}
\begin{proof}  Suppose that $r = p^k$ and  $\Delta_K(\zeta_r^i) = 0$.  Then the $r$--cyclotomic
polynomial, $\phi_r(t) = (t^{p^k} - 1)/ (t^{p^{k-1}} - 1) $ would divide $\Delta_K(t)$.  But
$\phi_r(1) = p$   while $\Delta_K(1) =\pm 1$.
\end{proof}

We now proceed with the proof of Theorem \ref{polytheorem}.

\begin{proof}[Proof of Theorem \ref{polytheorem}]

 According to Riley \cite{ri}  the order of the homology of the
$k$--fold cyclic branched cover  of a knot $K$ grows exponentially as a function of $k$ if the
Alexander polynomial has a root that is not a root of unity.  Hence, we only need to consider
the case that all irreducible factors of the Alexander polynomial are cyclotomic
polynomials,
$\phi_n(t)$.  Using Theorem \ref{polyproduct}, the result is reduced to the case that that
$\Delta_K(t) = \phi_n(t)$.  As in the proof of Corollary  \ref{rat}, $n$ cannot be a prime power.

An elementary argument using the resultant of polynomials (see for instance \cite{la})  gives
$$\prod_{i=0 }^{ p^k -1}
\phi_n(\zeta_{p^k} ^i) = \prod
 ((\omega_n  )^{p^k} -1)$$
where the second product is taken over all primitive $n$--roots of unity.  Let $g =
\mbox{gcd}(n,p^k)$ and let $m = n/g$.  One has that $\omega_n^{p^k} = \omega_m$ for some primitive
$m$--root of unity and with a bit of care one sees that the product can be rewritten as 
$$  \prod 
 ( \omega_m    -1)    ^b$$ where now the product is over all primitive $m$--roots of unity and
$b \ge 1$.  (Though we don't need it, a close examination shows that if $k$ is greater than or equal
to the maximal power of
$p$ in
$n$ then
$b = p^k - p^{k-1}$, otherwise $b= p^k$.)

If $n$ has three distinct prime factors then $m$ has at least two distinct prime factors and this
product is 1 (see for instance \cite[page 73]{la2}).  On the other hand, if $n$ has two distinct
prime factors, then by letting $p$ be one of those factors and letting $k$ be large, it is arranged
that
$m$ is a prime power and the product yields that prime and in particular is greater than 1.  This
concludes the proof of the first statement of  Theorem \ref{polytheorem}.

Finally, suppose that all cyclic branched covers of $K$ are homology spheres.  By the above discussion
we just need to show that no factor of the Alexander polynomial is $\phi_n(t)$ for any $n$.  But  from
Theorem \ref{polyproduct}  we see that if $\phi_n(t)$ divides the Alexander polynomial then the
$n$--fold cyclic branched cover would have infinite homology.  This concludes the proof.
\end{proof}

%%%%%%%%%%%%%%%%%%%%%%%%%%%%%%%%%%%%%%%%%%%%%%%%%%%%%%%%%%%%%%%%%%%%%%%%

\end{document}